# First-Order Intuitionistic Logic with Decidable Propositional Atoms

Alexander Sakharov
alex@sakharov.net
http://alex.sakharov.net

**Abstract**
First-order intuitionistic logic extended with the assumption about decidability of all propositional atoms combines classical properties in its propositional part and intuitionistic properties for derivable formulas not containing propositional symbols. Sequent calculus is used as a framework for investigating this extension. Constrained Kripke structures are introduced for modeling intuitionistic logic with decidable propositional atoms. The extent of the disjunction and existence properties is investigated.

## Introduction

The difference between intuitionistic and classical logic can be characterized by the Aristotelian law of excluded middle (LEM): $A \lor \neg A$. From the intuitionistic point of view, LEM is not acceptable because there is no general method for deciding whether $A(x)$ is true or false for arbitrary *x* if x ranges over an infinite domain. From the classical point of view, LEM is a valid and very important principle used everywhere in mathematics. As Hilbert stated [Hi]: "Taking the principle of excluded middle from the mathematician would be the same, proscribing the telescope to the astronomer or to the boxer the use of his fists. To prohibit existence statements and the principle of excluded middle is tantamount to relinquishing the science of mathematics altogether."

Which of the two logics - classical and intuitionistic - is better? This rhetorical question has been around for about a century. The matter of fact is that each has certain advantages over the other in certain dimensions. Both intuitionistic and classical logic have issues, too. Classical logic has decent models defined with using truth tables, but existence proofs are not constructive in it. Intuitionistic logic has the disjunction and existence properties and thus constructive proofs but the propositional fragment of intuitionistic logic while being finite does not have precise finite models.

Is it possible to combine the best of both worlds? Research in this domain was initiated by Godel and Tarski almost as soon as intuitionistic logic emerged. First of all, classical logic can be interpreted within intuitionistic one via Godel-Gentzen-Kolmogorov's double negation translation. Since derivable intuitionistic formulas constitute a subset of derivable classical formulas, the focus of this research has been on investigating properties of the logics lying between the two. Such logics are called intermediate or superintuitionistic. Intermediate logics are usually defined by adding one or more axiom schemas weaker than LEM to intuitionistic logic. For example, the Dummett logic is defined by adding the axiom schema: $(A \supset B) \lor (B \supset A)$. See [Fi] for an overview of intermediate logics. Many intermediate logics have the disjunction property [FM, Ma].

Recently, research in the area of combining features of classical and intuitionistic logic has shifted towards logics containing two different variants of connectives. One of them is intuitionistic, and the other is classical [Kr] Fibring logics is the most noticeable technique in this research [Ga]. This approach is less intuitive though. First, it is generally not clear where to apply which of the two paired connectives. Second, joining intuitionistic and classical axiomatizations with careless addition of some interaction axioms may make the combination logic collapse into classical logic [FH, SRC]. Precautions are necessary to avoid the collapsing problem.

Let us have a fresh look at this problem of combining features of intuitionistic and classical logic and try to do it staying within the framework of intuitionistic logic. Note that Brouwer [He] criticized LEM and thus classical logic only because LEM was abstracted from finite situations and extended without justification to statements about infinite collections. As Weyl [We] states in reference to Brouwer: "classical logic was abstracted from the mathematics of finite sets and their subsets". Propositional logic can be considered a part of the mathematics of finite sets because of availability of finite models using truth tables. Thus, LEM for propositional formulas is not really a target of intuitionistic criticism of classical logic.

The classical assumption that every propositional atom (propositional symbol, prime propositional formula) is either true or false but not both is so natural. This assumption translates into decidability of propositional atoms in terms of proof theory, and it leads to LEM for propositional formulas. Let us look what happens to first-order intuitionistic logic when assuming that all propositional atoms are decidable. Extending this assumption onto all prime predicate formulas is not quite natural though as discovered by intuitionists. Nevertheless, Heyting arithmetic is an intuitionistic theory in which basic arithmetic formulas are decidable.

First-order intuitionistic logic with decidable propositional atoms has a balanced mix of classical and intuitionistic characteristics. It coincides with classical logic in its propositional part, and it coincides with intuitionistic logic on the set of formulas either not containing propositional symbols or having occurrences of the same polarity for every propositional symbol. The disjunction property holds for disjunctions whose disjuncts do not share any propositional symbols. Some classes of formulas for which the existence property holds are also identified. Additionally, a weaker form of the existence property holds in intuitionistic logic with decidable propositional atoms.

Sequent calculus is used as a framework for studying first-order intuitionistic logic with decidable propositional atoms. Intuitionistic logic with decidable propositional atoms is sound and complete with respect to constrained Kripke structures introduced here. Absence of paired connectives makes this approach to combining classical and intuitionistic features more natural and intuitive. This approach does not suffer the collapsing problem either. In contrast to intermediate logics, there is no introduction of non-intuitionistic axioms applicable to prime predicate formulas that are relevant to infinite collections.

First-order intuitionistic logic with decidable propositional atoms is supposed to be a theory somewhat acceptable by both camps: classical and intuitionistic. On one hand, it is fully classical in its propositional part. On the other hand, LEM, double negation elimination and their variants are not extended onto statements about infinite collections.

Note that research with slightly similar motivations and goals was recently done by H. Kurokawa [Ku]. He investigated propositional logic with two types of atoms: classical and intuitionistic.

And LEM is applicable to classical atoms only. Although, Kurokawa's paper does not cover first-order issues. From the perspective of Kurokawa's research, first-order logic with classical propositional atoms and intuitionistic predicate atoms is investigated here.

## Background and Definitions

Investigation of properties of first-order intuitionistic logic with decidable propositional atoms will be based on using sequent calculus as a framework. We use the phrase 'predicate symbol' only in application to symbols whose arity is more than zero because symbols of arity zero are decidable and thus have very different properties than the rest of formulas. We will use capital Latin letters to denote formulas whereas capital Greek letters will denote (possibly empty) sequences of formulas.

**Definition**. A formula is called propositional if it contains only propositional symbols. A formula is called purely predicate if all its symbols are predicate symbols. A formula is called mixed if it contains both propositional and predicate symbols.

Let us repeat the Gentzen's formulation of intuitionistic logic (LJ) from [Ge]. A sequent is an expression $\Gamma \vdash \Lambda$ where $\Gamma$ is a (possibly empty) sequence of formulas and $\Lambda$ is either a formula or an empty sequence. $\Gamma$ is called antecedent, $\Lambda$ is called succedent. The only axiom schema of LJ is:
$A \vdash A$

The rules of inference are: thinning, contraction, exchange, cut, pairs of rules for negation, implication, existential quantifier, universal quantifier, and triples of rules for conjunction and disjunction. The pairs and triples include rules operating on both antecedents and succedents. These rules are.

Thinning:
$$\frac{\Gamma \vdash \Lambda}{A, \Gamma \vdash \Lambda} \qquad \frac{\Gamma \vdash}{\Gamma \vdash A}$$

Contraction:
$$\frac{A, A, \Gamma \vdash \Lambda}{A, \Gamma \vdash \Lambda}$$

Exchange:
$$\frac{\Pi, A, B, \Gamma \vdash \Lambda}{\Pi, B, A, \Gamma \vdash \Lambda}$$

Cut:
$$\frac{\Gamma \vdash C \qquad C, \Pi \vdash \Lambda}{\Gamma, \Pi \vdash \Lambda}$$

&:
$$\frac{\Gamma \vdash A \qquad \Gamma \vdash B}{\Gamma \vdash A\&B} \qquad \frac{A, \Gamma \vdash \Lambda}{A\&B, \Gamma \vdash \Lambda} \qquad \frac{B, \Gamma \vdash \Lambda}{A\&B, \Gamma \vdash \Lambda}$$

∨:
$$\frac{A, \Gamma \vdash \Lambda \quad B, \Gamma \vdash \Lambda}{A \vee B, \Gamma \vdash \Lambda} \qquad \frac{\Gamma \vdash A}{\Gamma \vdash A \vee B} \qquad \frac{\Gamma \vdash B}{\Gamma \vdash A \vee B}$$

¬:
$$\frac{A, \Gamma \vdash}{\Gamma \vdash \neg A} \qquad \frac{\Gamma \vdash A}{\neg A, \Gamma \vdash}$$

⊃:
$$\frac{A, \Gamma \vdash B}{\Gamma \vdash A \supset B} \qquad \frac{\Gamma \vdash A \quad B, \Pi \vdash \Lambda}{A \supset B, \Gamma, \Pi \vdash \Lambda}$$

∀:
$$\frac{\Gamma \vdash F(a)}{\Gamma \vdash \forall x\, F(x)} \qquad \frac{F(a), \Gamma \vdash \Lambda}{\forall x\, F(x), \Gamma \vdash \Lambda}$$

∃:
$$\frac{\Gamma \vdash F(a)}{\Gamma \vdash \exists x\, F(x)} \qquad \frac{F(a), \Gamma \vdash \Lambda}{\exists x\, F(x), \Gamma \vdash \Lambda}$$

Here variable $a$ occurs free in $F$, $F(x)$ is obtained from $F(a)$ by replacing occurrences of $a$ by $x$. Variable $a$ in the ∀-succedent rule and in the ∃-antecedent rule is called eigenvariable. It must not occur in the lower sequents of the respective rules.

Typically, two constants - ⊤ and ⊥ - are introduced as well. ⊤ denotes a 'valid' formula whereas ⊥ denotes an 'unsatisfiable' formula. Formally, these constants are introduced by the following axioms:

$\vdash \top \qquad \bot \vdash$

Formula $F$ is called decidable if $F \vee \neg F$ holds in a respective theory. From the proof-theoretic point of view, assumption that all propositional atoms are decidable means addition of axioms $A \vee \neg A$ for all propositional symbols $A$. Basically, this assumption amounts to extension of LJ with the axiom schema

$\vdash A \vee \neg A$

that holds for all propositional atoms. Now, let us introduce an alternative definition for first-order intuitionistic logic augmented with the decidability assumption for propositional atoms. This alternative definition is represented as LJ extended with an inference rule. This extension is called LJ+ and will be extensively used later.

**Definition.** LJ+ is LJ extended with one additional rule:
$$\frac{N, \Gamma \vdash \Lambda \qquad \neg N, \Gamma \vdash \Lambda}{\Gamma \vdash \Lambda}$$
It will be called neutralization rule. Only propositional formulas $N$ are allowed in the neutralization rule.

LJ+ is equivalent to intuitionistic logic extended with decidable propositional atoms, that is, the sets of derivable sequents are the same. First, consider an arbitrary LJ+ derivation. As well-known, decidability of propositional atoms implies LEM for all propositional formulas, i.e. N $\lor \neg$N is derivable for any propositional formula N. Applications of the neutralization rule are eliminated as follows:

```
                    N, Γ |-- Λ        ¬N, Γ |-- Λ
                    ----------------------------------
|-- N∨¬N                    N∨¬N, Γ |-- Λ
-----------------------------------------------
            Γ |-- Λ
```

The remaining derivation parts are unaffected. Second, consider an arbitrary intuitionistic derivation including instances of the aforementioned axiom schema. Sequents |-- N$\lor\neg$N, where N is a propositional symbol, can be replaced by the following derivations using the neutralization rule:

```
N |-- N          ¬N |-- ¬N
------------     --------------
N |-- N∨¬N       ¬N |-- N∨¬N
------------------------------------
            |-- N∨¬N
```

Apparently, every formula derivable in LJ+ is also derivable in LK. LK is Gentzen's formulation of classical logic. Gentzen [Ge] alternatively defines LK as LJ (i.e. all sequents have no more than one formula in the succedent) plus axiom schema |-- A$\lor\neg$A that is applicable to any formula A. LJ+ is consistent because the set of derivable formulas of LJ+ is a subset of those of LK, and LK is known to be consistent.

There exist formulas provable in LJ+ but not provable in LJ. R$\lor\neg$R for propositional symbol R is an example of such formula holding LJ+. See [Ge] for the proof that R$\lor\neg$R does not hold in LJ. Also, there exist formulas provable in LK but not in LJ+. This will be proven later - see section 'Intuitionistic Features'. $\exists x\, P(x) \lor \neg \exists x\, P(x)$ is an example of a classically provable formula that does not hold in LJ+.

Relationship between intuitionistic logic with decidable propositional atoms and pure intuitionistic logic can be characterized as the following. Formula F is derivable in LJ+ if and only if, for some classically derivable propositional formula A, A$\supset$F is derivable in LJ. The 'only if' part is a corollary of the deduction theorem for intuitionistic logic and the fact that B$\supset$(C$\supset$D) is equivalent to (B&C)$\supset$D in intuitionistic logic. Formula A is the conjunction of all axioms N$\lor\neg$N used in the derivation of F. As shown later in the section 'Classical Features', the propositional fragment of LJ+ is equivalent to classical propositional logic. As for the 'if' part, if A is classically derivable propositional formula, it is also derivable in LJ+. LJ+ is closed under modus ponens because a well-know proof for sequent calculus works for LJ+ too. Closeness of LJ+ under modus ponens guarantees that F is also derivable in LJ+.

## Sample Derivations

This section illustrates what kind of derivations can be conducted in intuitionistic logic with decidable propositional atoms. Three simple examples are presented below. The first example uses decidability directly whereas the other two are derivations containing the neutralization rule. More sample derivations can be found in sections 'Classical Features' and 'Disjunction and

Existence Properties'. Note that as a part of the second example, double negation elimination is derived for propositional formulas.

**Example 1**

```
  |-- A\/ ¬A                         …
---------------------      ----------------------
∀x B(x) |-- A\/ ¬A         ∀x B(x) |-- ∃x B(x)
----------------------------------------------------
∀x B(x)  |-- (A \/ ¬A) & ∃x B(x)
------------------------------------
|--- ∀x B(x) ⊃ ((A \/ ¬A) & ∃x B(x))
```

**Example 2**

```
                 ¬A |-- ¬A
                 ------------
A |-- A          ¬A, ¬¬A |--
--------         ----------------
A, ¬¬A |-- A    ¬A, ¬¬A |-- A
----------------------------------------
        ¬¬A |-- A                  ∀x B(x)  |-- ∀x B(x)
        -----------                ---------------------------
        ¬¬A |-- A \/ ∀x B(x)   ∀x B(x) |-- A \/ ∀x B(x)
        -----------------------------------------------------
            ¬¬A \/ ∀x B(x)  |-- A \/ ∀x B(x)
            -----------------------------------------
            |-- (¬¬A \/ ∀x B(x)) ⊃ (A \/ ∀x B(x))
```

**Example 3**

```
A |-- A                                        ¬A |-- ¬A
--------                                       -------------
A |-- A \/ ∀x P(x)                             ¬A |-- ¬A \/ ∀x Q(x)
--------------------                           ------------------------
A |-- (A \/ ∀x P(x))\/(¬A \/ ∀x Q(x))    ¬A |-- (A \/ ∀x P(x))\/(¬A \/ ∀x Q(x))
-----------------------------------------------------------------------------------
        |-- (A \/ ∀x P(x))\/(¬A \/ ∀x Q(x))
```

## Classical Features

The propositional fragment of LJ+ is equivalent to classical propositional logic, that is, to the propositional fragment of LK. As mentioned earlier, LK can be alternatively defined as LJ plus axiom schema |-- A\/¬A. Given that formulas A\/¬A are derivable in LJ+ for all propositional formulas, it is clear that every LK derivation (by the above definition) containing only propositional symbols is also a valid LJ+ derivation.

There exist mixed formulas provable in LJ+ but not provable in LJ. Consider formula (R\/¬R) \/ P(x) in which both R and P are atomic. This formula is apparently derivable in LJ+ because R \/¬R is derivable in LJ+. Assume it is derivable in LJ. Due to the disjunction property, either R \/¬R or P(x) should be derivable in LJ. P(x) could not be derivable in LJ because it could not be the bottom sequent of any rule except thinning of $\vdash$ which in turn is not derivable in LJ. It is known that $R \lor \neg R$ is not derivable in LJ either.

The following theorem is characteristic of classical logic, and it does not hold in intuitionistic logic. Interestingly, this theorem holds in LJ+. If R is a propositional symbol occurring in formula A, then let $A\{R|F\}$ denote the formula obtained from A by replacing all occurrences of R with formula F.

**Theorem**. If R is a propositional symbol occurring in formula A, then A is derivable in LJ+ if and only if both $A\{R|\top\}$ and $A\{R|\bot\}$ are derivable in LJ+.

**Proof**. The 'only if' part is proved by replacing all occurrences of R by $\top$ or by $\bot$ in A's derivation, respectively. Apparently, the derivation remains intact. Now consider derivations of $A\{R|\top\}$ and $A\{R|\bot\}$ and the occurrences of $\top$ and $\bot$ that replaced R. The constants $\top$ and $\bot$ introduced in derivations of $A\{R|\top\}$ and $A\{R|\bot\}$ by other means than the axioms
$\vdash \top$    and    $\bot \vdash$
can be replaced by R retaining the derivation intact because the only other means are thinning, axioms $\top \vdash \top$, $\bot \vdash \bot$. If $\top$ is introduced by axiom $\vdash \top$ in derivation of $A\{R|\top\}$, let us replace this axiom by axiom
$R \vdash R$
If $\bot$ is introduced by axiom $\bot \vdash$ in derivation of $A\{R|\bot\}$, let us replace this axiom by
$R, \neg R \vdash$
which is derivable from the above axiom in one step.

After that, let us transform derivations of $A\{R|\top\}$ and $A\{R|\bot\}$ so that R ($\neg R$) is added as the rightmost formula to the antecedents of all sequents below the replaced axioms for $\top$ ($\bot$) and to the antecedent of their counterpart sequents in the upper parts of rules. This transformation is done recursively top down to the endsequent. If a rule has one upper sequent and R ($\neg R$) has been already added as the rightmost formula to the antecedent of its upper sequent, then let us add R ($\neg R$) as the rightmost formula in the antecedent of the lower sequent. Apparently, this rule remains intact. If a rule has two upper sequents and R ($\neg R$) has been already added as the rightmost formula to the antecedent of one of its upper sequents, then let us add R ($\neg R$) to the counterpart upper sequent by injection of thinning and exchanges in order to make R ($\neg R$) the rightmost formula of the antecedent of this sequent. If we also add R ($\neg R$) as the rightmost formula to the antecedent of the lower sequent, then this rule remains intact. In case R ($\neg R$) has already been added to both upper sequents, let us add R ($\neg R$) as the rightmost formula to the antecedent of the lower sequent, and this rule remains intact as well. As a result of this transformation, both derivations under consideration remain correct. The endsequent in the derivation of $A\{R|\top\}$ will have one formula R in the antecedent whereas the endsequent in the derivation of $A\{R|\bot\}$ will have one formula $\neg R$ in the antecedent. These two derivations are merged by application of the neutralization rule to R and $\neg R$. The resulting endsequent is A. ♦

One other interesting classical feature of LJ+ is that $N \supset F$ is equivalent to $\neg N \lor F$ for any formula F and any propositional formula N. This is proved directly by constructing simple derivations for the implications in both directions. In the following two figures presenting these derivations, double lines denote application of more than one rule.

```
N |-- N     F |-- F
------------------
N⊃F, N |-- F                    ¬N |-- ¬N
========                        =======
N, N⊃F |-- ¬N∨F                 ¬N, N⊃F |-- ¬N∨F
-----------------------------------------------------------
        N⊃F |-- ¬N∨F

N |-- N                 F |-- F
=====                   ====
¬N, N |-- F             F, N |-- F
---------------------------------------
        ¬N∨F, N |-- F
        ----------------
        ¬N∨F |-- N⊃F
```

## Weak Cut Elimination

The neutralization rule is another discharge rule along with the cut rule. The cut rule can be eliminated from LJ+ derivations but the neutralization rule cannot. Therefore, this weak form of cut elimination achievable in LJ+ does not imply the subformula property. The subformula property holds in application to predicate formulas only. Nevertheless, admissibility of cut in LJ+ will be used later.

**Theorem**. Every LJ+ derivation can be transformed into another derivation with the same endsequent and in which the cut rule does not occur.

**Proof**. We use Gentzen's cut elimination proof for LJ [Ge] and augment it with a few additional transformations for derivation fragments involving the neutralization rule. Transformation of derivations with cut into derivations with mix remains unchanged. The mix rule is the following:

```
Γ |-- Θ        Π |-- Λ
---------------------------
     Γ, Π* |-- Θ*, Λ
```

Both Θ and Π contain so-called mix formula A. Θ* and Π* are obtained from Θ and Π by removing A from the respective formula sequences.

Apparently, formula N from the neutralization rule

```
N, Γ |-- Λ        ¬N, Γ |-- Λ
---------------------------------
            Γ |-- Λ
```

cannot be the mix formula of mix immediately below this rule. Therefore, Gentzen's proof for the case when both the left and the right rank equal one does not require any augmentation because the upper sequents of mix cannot be lower sequents of the neutralization rule.

Let us look at the part of Gentzen's proof dealing with lowering either the left or right rank. In addition to Gentzen's proof, we only have to consider two cases in which the left or right upper sequent of mix is the lower sequent of the neutralization rule. If the right upper sequent of mix is

the lower sequent of the neutralization rule

```
                    N, Γ |-- Λ        ¬N, Γ |-- Λ
                    ---------------------------------- neutralization
         Π |-- A              Γ |-- Λ
         ---------------------------------------- mix
                   Π, Γ* |-- B
```

then it is transformed into

```
         Π |-- A   N, Γ |-- Λ           Π |-- A   ¬N, Γ |-- Λ
         ---------------------- mix     ------------------------ mix
         Π, N*, Γ* |-- Λ                 Π, ¬N*, Γ*  |-- Λ
         ------------------              ----------------------
               thinnings and exchanges
         ------------------              --------------------
         N, Π, Γ* |-- Λ                  ¬N, Π, Γ*  |-- Λ
         ------------------------------------------------------------- neutralization
                   Π, Γ* |-- Λ
```

Γ* is obtained from Γ by deleting all occurrences of A. N* is N if A is different from N or nil otherwise. ¬N* is ¬N if A is different from ¬N or nil otherwise. In this transformed derivation, the rank of each of the two newly introduced mix rules is smaller than that of the original mix rule.

If the left upper sequent of mix is the lower sequent of the neutralization rule

```
         N, Γ |-- B        ¬N, Γ |-- B
         -------------------------------- neutralization
                 Γ |-- B              Π |-- Λ
                 ------------------------------------------ mix
                         Γ, Π* |-- Λ
```

then it is transformed into

```
         N, Γ |-- B      Π |-- Λ          ¬N, Γ |-- B      Π |-- Λ
         ------------------------------ mix  ------------------------------ mix
                N, Γ, Π* |-- Λ                     ¬N, Γ, Π*  |-- Λ
                ------------------------------------------------------------- neutralization
                                   Γ, Π* |-- Λ
```

Π* is obtained from Π by deleting all occurrences of B. In this transformed derivation, the rank of each of the two newly introduced mix rules is smaller than that of the original mix rule. ♦

Seldin [Se] achieves a similar result by detouring via natural deduction and normalizing proofs there. This weak form of cut elimination does not imply the subformula property because propositional formulas can be discarded by application of the neutralization rule. The subformula property holds for predicate formulas though. Namely, every predicate formula in a cut-free derivation is a subformula of some formula in the endsequent. The proof is similar to the proof of the subformula property given by Gentzen in [Ge]. Just note that none of predicate formulas can be discarded by application of the neutralization rule.

Derivations in LJ+ without cut have another interesting property. If all formulas in a rule are propositional, then all formulas above this rule in a cut-free derivation are propositional. This follows from the fact that only propositional formulas can be discarded in derivations without cut.

# Constrained Kripke Structures

Let L denote a first-order language. A Kripke structure (frame) K consists of a partially ordered set K of nodes and a function D assigning to each node k in K a domain D(k) such that if k≤k', then D(k) ⊂ D(k'). In addition, a valuation T(k) maps k to a set of atomic expressions of L over D(k) for all k∈ K. T satisfies the following: if k ≤ k' then T(k) ⊂ T(k').

A forcing relation is defined recursively for L as follows:

k forces atomic expression $A(d_1,…,d_n)$ iff $A(d_1,…,d_n) \in T(k)$
k forces A&B if k forces A and k forces B
k forces A⋁B if k forces A or k forces B
k forces A⊃B if, for every such k' that k≤k', if k' forces A then k' forces B
k forces ¬A if for no such k' that k≤k', does k' force A
k forces ∀x A(x) if for every such k' that k≤k' and every d ∈ D(k'), k' forces A(d)
k forces ∃x A(x) if for some d ∈ D(k), k forces A(d)
⊤ is forced by all nodes in all Kripke structures
⊥ is not forced by any node in any Kripke structure

Any such forcing relation is known to be consistent and monotone. Intuitionistic logic is sound and complete with respect to Kripke structures [Mi, TD]
**Consistency**: For no sentence A and no k does k force both A and ¬A.
**Monotonicity**: If k ≤ k' and k forces A then k' forces A.
**Soundness**: If formula A is derivable in LJ, then A is forced at all nodes in every Kripke structure.
**Completeness**: If A is forced at all nodes in every Kripke structure, then A is derivable in LJ.

For modeling LJ+, we consider a subset of Kripke structures defined by the following constraint: if k ≤ k' and Q is a propositional symbol, then either both T(k) and T(k') contain Q or they both do not contain Q. Let us call this subset constrained Kripke structures. Our goal is to select such Kripke structures that behave like the classical truth tables on the propositional fragment of LJ+. The following theorem justifies the choice of constrained Kripke structures.

**Theorem**. For any propositional formula A and any node k of any constrained Kripke structure K, either k forces A or k forces ¬A. For any two nodes k and k' of K, if k ≤ k' and k' forces A, then k forces A.
**Corollary.** For any node k of any constrained Kripke structure K and any propositional formulas A and B, k forces A⊃B if and only if 'k forces A' implies 'k forces B', k forces ¬A if and only if k does not force A.

**Proof**. Both statements of this theorem are proven by induction on the size of propositional formula in question. To start proving the first statement of the theorem, consider the case that A is atomic, i.e. it is propositional symbol Q. If T(k) is contains Q, then k forces A. If T(k) does not contain A, then k forces ¬A because T(k') does not contain Q for all such k' that k≤k'. Suppose A is B&C, and the first statement holds for B and C. For any k, it sufficient to consider four cases: k forces A and B, k forces A and ¬B, k forces ¬A and B, k forces ¬A and ¬B. If k forces A and B, then it also forces A&B. In the three remaining cases, for no k' such that k≤k', does k' force both A and B. Therefore, k forces ¬(A&B). Suppose A is B⊃C, and the first statement holds for B and C. Again, it sufficient to consider four cases: k forces A and B, k forces A and ¬B, k forces ¬A and B, k forces ¬A and ¬B. If k forces A and B, or k forces ¬A and B, or k forces ¬A and ¬B, then k forces B⊃C due to the monotonicity of Kripke structures. If k forces A and ¬B, then for no k' such that k≤k', does k' force A⊃B. Therefore, k forces ¬(A⊃B). Cases of other connectives are considered similarly.

To prove the second statement of the theorem, again start from considering the case that A is atomic, i.e. it is propositional symbol Q. If k≤ k' and k' forces Q, then k forces Q by definition of the constrained Kripke structures. Suppose A is B$\lor$C and the second statement holds for B and C. If k' forces B$\lor$C, then k' forces B or it forces C. If k' forces B and k≤k', then k forces B by the induction assumption. Therefore, k forces B$\lor$C as well. Considering the case that k' forces C leads to the same conclusion. Now suppose A is ¬B and the second statement holds for B. Let k ≤k' and k' forces ¬B. Consider arbitrary k" such that k≤k". If k" forces B, then k forces B by the induction assumption. By monotonicity, k' forces B, which contradicts the fact that k' forces ¬B. Therefore, k" does not force B and thence k forces ¬B. Cases of other connectives are considered similarly. ♦

**Definition.** Sequent Γ |-- Λ is called valid if for every node k in every constrained Kripke structure K, either there is such formula A from Γ that k does not force A or k forces Λ provided that Λ is not empty.

**Theorem** (Soundness of LJ+ with respect to constrained Kripke structures). Any sequent Γ |-- Λ derivable in LJ+ is valid.
**Corollary**. If formula A is derivable in LJ+, then A is forced at all nodes in every constrained Kripke structure.

**Proof.** Let us look at derivations in LJ+ without cut. We use induction on the height of derivation of Γ |-- A. Basically, we only have to check out that the neutralization rule preserves validity, that is, the lower sequent of the neutralization rule is valid whenever the upper sequents are. It is well known that all other rules that may appear in derivations preserve validity [Mi], and apparently all instances of the axiom schema A |-- A are valid as well. Therefore, induction step can be made in application to all rules except neutralization.

Let us look at the neutralization rule:
N, Γ |-- Λ        ¬N, Γ |-- Λ
---------------------------------
          Γ |-- Λ

As established earlier, for any node k of any constrained Kripke structure K, either k forces N or k forces ¬N. If k forces N, consider N, Γ |-- Λ, otherwise consider ¬N, Γ |-- Λ. By the induction assumption applied to N, Γ |-- Λ (or to ¬N, Γ |-- Λ, respectively), either there is such formula A from N, Γ (¬N, Γ ) that k does not force A or k forces Λ provided that Λ is not empty. A is from Γ because it cannot be N (¬N). Therefore, Γ |-- Λ is valid. ♦

Completeness of LJ+ with respect to constrained Kripke structures will be proven later in the section 'Intuitionistic Features'. Soundness of LJ+ with respect to constrained Kripke structures is a useful tool in showing that certain formulas cannot be proved in LJ+. For instance, let us prove that ∃x P(x) $\lor$ ¬∃x P(x) does not hold in LJ+. Consider a constrained Kripke structure with two elements k≤ m, D(k) = D(m) = { 0 }, T(k) is empty and T(m) consists of one element P(0). ∃x P(x) is not forced at k because T(k) is empty. ¬∃x P(x) is not forced at k because ∃x P(x) is forced at m and k≤m. Therefore, ∃x P(x) ∨ ¬∃x P(x) is not derivable in LJ+ because k does not force it.

## Intuitionistic Features

Once we established the soundness of LJ+ with respect to constrained Kripke structures, we can address the issue of intuitionistic features of LJ+.

**Lemma**. The purely predicate fragment of LJ+ is equivalent to LJ, that is, every purely predicate formula derivable in LJ+ is also derivable in LJ.

**Proof.** Consider purely predicate formula F derivable in LJ+. Consider an arbitrary Kripke structure K (not necessarily constrained) and any node k∈K. Let us build another Kripke structure K' whose only difference from K is that all propositional symbols are removed from valuations T(k) for all k∈K. K' satisfies the definition of constrained Kripke structures simply because valuations do not have propositional symbols. Since F is purely predicate, k forces F in K if and only if k forces F in K' for any k∈K(K'). By soundness of LJ+ with respect to constrained Kripke structures, k forces F for all k∈K'. Therefore, k forces F in K as well. By completeness of LJ with respect to Kripke structures, F is derivable in LJ. ♦

**Definition**. Positive (negative) subformula occurrences in propositional formula A are defined recursively by the following rules:
- A itself is positive.
- If subformula B∨C or B&C is positive (negative), then A and B are positive (negative).
- If subformula ¬B is positive (negative), then B is negative (positive).
- If subformula B⊃C is positive (negative), then B is negative (positive) and C is positive (negative).
- If subformula ∀x B(x) or ∃x B(x) is positive (negative), then B(x) is positive (negative).

**Definition.** A formula is called unipolar if occurrences of every propositional symbol in it are either all positive or all negative.

The next theorem shows that LJ+ coincides with LJ on the class of unipolar formulas. Given that the simplest formulas having both negative and positive occurrences of propositional symbols are LEM and double negation elimination, this class perhaps gives a reasonable syntactically-defined approximation of the boundaries within which LJ+ and LJ coincide. Suppose B is a subformula of A. Consider a particular occurrence of B in A. Let $A_{B|C}$ denote the formula obtained from A by replacing the given occurrence of B by C.

**Theorem**. Every unipolar formula derivable in LJ+ is also derivable in LJ. There exist formulas derivable in LK but not derivable in LJ+.
**Corollary**. The disjunction and existence properties hold in LJ+ in application to unipolar formulas.

**Proof.** The proof of this theorem is by induction on the number of propositional symbols in formula F. Base case: F has no propositional symbols. The previous lemma gives a proof of the base case. Induction step: suppose this theorem holds for formulas having not more than n propositional symbols, and F has n+1 propositional symbols. Suppose R is one of F's propositional symbols, all occurrences of R are positive. The case of negative occurrences is similar to this one. If F is derivable in LJ+, so is F{R|⊥}. By the induction hypothesis, if F{R|⊥} is derivable in LJ+, then it is derivable in LJ because it contains n propositional symbols.

LJ has the following well-known properties: if $A \supset B$ is derivable and A is a positive subformula of F, then $F \supset F_{A|B}$ is derivable; if $C \supset A$ is derivable and A is a negative subformula of F, then $F \supset F_{A|C}$ is derivable. Consider such chain of formulas that the first formula is F{R|⊥}, every next formula in this chain is received by replacing one occurrence of ⊥ back by R. The occurrences of ⊥ that have not been generated by the original replacement of R are not involved. The last formula in this chain is F. Since $\bot \supset R$, every formula in this chain except the last one implies the next formula. Since implication is transitive, $F\{R|\bot\} \supset F$ is derivable in LJ and so is F due to closeness under modus ponens.

Now let us prove the second statement of this theorem. Consider formula $(\exists x\, P(x) \lor \neg \exists x\, P(x)) \lor R$ where R is a propositional symbol. This formula is provable in LK because $\exists x\, P(x) \lor \neg \exists x\, P(x)$ is. Look at the constrained Kripke structure presented earlier. It has two elements $k \leq m$, $D(k) = D(m) = \{0\}$, T(k) is empty and T(m) consists of one element P(0). $\exists x\, P(x)$ is not forced at k because T(k) is empty. $\neg \exists x\, P(x)$ cannot not forced at k because $\exists x\, P(x)$ is forced at m and $k \leq m$. R is not forced at k either because T(k) does not contain R. Thus $(\exists x\, P(x) \lor \neg \exists x\, P(x)) \lor R$ is not forced at k. The soundness of LJ+ implies that $(\exists x\, P(x) \lor \neg \exists x\, P(x)) \lor R$ is not derivable. ♦

This theorem gives an indirect proof of admissibility of the neutralization rule in derivations of unipolar formulas. It follows from admissibility of cut in intuitionistic logic. The neutralization rule is not admissible in LJ+ though. If it had been admissible, then all derivations would have been purely intuitionistic.

**Definition**. Occurrence of formula B in formula A is called strictly positive if there are no outer connectives for B except for & and $\lor$ and this occurrence of B is not in the scope of any quantifier in A.

**Theorem**. If A is classically provable purely predicate formula that is not derivable in LJ+, B is strictly positive subformula occurrence in A, and C is a propositional formula not derivable in LJ+, then $A_{B|B \lor C}$ does not hold in LJ+. If D is any propositional formula, then $A_{B|B\&D}$ does not hold in LJ+.

**Proof**. Since A is not derivable in LJ, there is Kripke structure K and node k in it such that k does not force A due to completeness of LJ with respect to Kripke structures. It is safe to assume that K does not contain propositional symbols (see the proof of the previous theorem). K is a constrained Kripke structure. If k forces B, then k may or may not force B&D. Given the rules for calculating forcing relationship for connectives & and $\lor$, k still does not force $A_{B|B\&D}$ in either case. If k does not force B, then it does not force B&D either. Again by forcing definition for & and $\lor$, k does not force $A_{B|B\&D}$. Therefore, $A_{B|B\&D}$ is not derivable in LJ+ because of soundness of LJ+ with respect to constrained Kripke structures.

Now consider the case with disjunction. Let us pick a classical logic model in which C is false. Extend K by adding elements to all D(k) and T(k) in accordance with the selected classical model, i.e. propositional symbol R is added to D(k) and T(k) for all $k \in K$ if and only if R is true in the selected classical model. This extended structure K' is apparently constrained. If k forces B in K', then k forces B$\lor$C too. If k does not force B, then k does not force B$\lor$C in K' because C is not forced by any node from K'. By forcing definition for & and $\lor$, k does not force $A_{B|B\lor C}$. This

formula is not derivable in LJ+ because of soundness of LJ+ with respect to constrained Kripke structures. ♦

This theorem assures that variants of LEM for predicate formulas do not hold in LJ+. If R is a propositional symbol, then neither $(\exists x\, P(x) \lor R) \lor \neg \exists x\, P(x)$ nor $(\exists x\, P(x)\, \&\, (R \supset R)) \lor \neg \exists x\, P(x)$ holds in LJ+. Likewise, formulas combining double negation elimination for predicate formulas with propositional symbols do not hold in LJ+ either. The following simple theorem gives a hint for constructing such formulas.

**Theorem**. If $A \supset B$ is a classically derivable purely predicate formula that is not derivable in LJ+, C is a propositional formula derivable in LJ+, D is a propositional formula not derivable in LJ+, then neither $A \& C \supset B$ nor $A \supset (B \lor D)$ holds in LJ+.

**Proof**. Since $A \supset B$ is not derivable in LJ, there is Kripke structure K and such node k in it that k does not force $A \supset B$ due to completeness of LJ with respect to Kripke structures. Let us select such $k \in K$ that k forces A and does not force B. As above, it is safe to assume that K does not contain propositional symbols, and thus, K is a constrained Kripke structure. Since C is derivable in LJ+, k forces A&C as well. Therefore, k does not force $A \& C \supset B$ and hence $A \& C \supset B$ is not derivable in LJ+ because of its soundness with respect to constrained Kripke structures. Using the same argument as in the previous theorem, we can extend K to another constrained Kripke structure in such a way that k does not force $B \lor D$. Therefore, k does not force $A \supset (B \lor D)$ and it is not derivable in LJ+. ♦

Here are a couple of sample variants of double negation elimination: $(\neg\neg(\exists x\, P(x))\, \&\, (A \supset A \lor B)) \supset \exists x\, P(x)$ and $\neg\neg \exists x\, P(x) \supset (\exists x\, P(x) \lor A)$. These formulas are not derivable in LJ+ as guaranteed by the last theorem.

**Theorem** (Completeness of LJ+ with respect to constrained Kripke structures). If formula F is forced at any node k in any constrained Kripke structure K, then F is derivable in LJ+.

**Proof**. The proof of this theorem is by induction on the number of propositional symbols in F. Base case: F has no propositional symbols. Consider an arbitrary Kripke structure K (not necessarily constrained) and any node $k \in K$. As explained earlier, it is safe to assume that K does not contain propositional symbols. Therefore, K is a constrained structure and k forces F. By the completeness of LJ with respect to Kripke structures, F is derivable in LJ and hence in LJ+.

Induction step: suppose this theorem holds for formulas having not more than n propositional symbols, and F has n+1 propositional symbols. If R is one of F's propositional symbols, then F is derivable in LJ+ if both $F\{R|T\}$ and $F\{R|\bot\}$ are derivable in LJ+ (by the theorem from the section 'Classical Features'.) Consider an arbitrary constrained Kripke structure K not containing R in any domain D, and consider any node $k \in K$. Since k forces F, k also forces $F\{R|\bot\}$ because forcing is identical for k paired up with any atomic subformula of F and k paired up with the respective atomic subformula of $F\{R|\bot\}$, and the structure of both formulas is the same. Now extend K with R belonging to all D(k) and all T(k). Since k forces F, k also forces $F\{R|T\}$ by the same argument. Apparently, presence of R in D(k) and T(k) does not affect forcing for $F\{R|T\}$ and $F\{R|\bot\}$ because they do not contain R. Therefore, $F\{R|T\}$ and $F\{R|\bot\}$ are forced by any k in any constrained Kripke structure, and thus they are derivable in LJ+ by the induction hypothesis, and so is F. ♦

## Disjunction and Existence Properties

The following disjunction and existence properties are perhaps the two most important features of intuitionistic logic:

**Disjunction:** If $A \lor B$ is derivable, then either A or B is derivable
**Existence:** If $\exists x\, A(x)$ is derivable, then $A(t)$ is derivable for some term t

The existence property is more interesting when logic formulations include functional symbols and object constants and therefore one can speak of terms. In our simplistic formulation without functional symbols and even without object constants, the existence property is reduced to a statement about a variable as opposed to a term.

Once the boundaries of LJ within LJ+ have been characterized, it is clear that both these properties hold for the formulas within the LJ boundaries. Apparently, the disjunction property does not hold for the entire LJ+ because it does not hold for classical propositional logic. The existence property does not hold in LJ+ either. Another question is: What is the scope of LJ+ outside LJ within which the disjunction and existence properties hold? The following two theorems give an insight on the extent of the disjunction and existence properties.

**Lemma.** If formula $F \lor G$ is derivable in LJ+, F and G do not have common propositional symbols, then either F or G is derivable in LJ+.

**Proof.** This proof is based on 'the gluing method' suggested by Kurokawa [Ku]. It is similar to Kurokawa's proof but applies to a different framework and yields a somewhat stronger result. Let F' denote the set of propositional symbols in F and G' denote the set of propositional symbols in G. If neither F nor G is derivable in LJ+, then there are such constrained Kripke structures K and M that F is not forced at some node $k' \in K$ and G is not forced at some node $m' \in M$. Let us 'glue' these two structures by adding a new node n such that $n \leq k$ for all $k \in K$, $k' \leq k$ and $n \leq m$ for all $m \in M$, $m' \leq m$. No non-propositional atomic expression $A(d_1, \ldots, d_n)$ is forced at n. Let the set of propositional symbols forced at n is the union of the propositional symbols forced at k' and the propositional symbols forced at m'. Also, let us add these propositional symbols forced at m' to the set of forced atomic expressions for every node in K, and let us add the propositional symbols forced at k' to the set of forced atomic expressions for every node in M.

This newly built combination of n, K and M is a constrained Kripke structure. Note that no pair $\langle k, m \rangle$ where $k \in K$ and $m \in M$ is ordered; $k \leq n$, $m \leq n$ do not hold either. Adding propositional symbols not occurring in F to the set of forced atomic expressions for all nodes of K does not affect valuation monotonicity. Neither this addition changes forcing relationships for nodes from K. The same is true about G and M, respectively. Similarly, all propositional symbols forced at n are also forced at nodes $k \in K$, $k' \leq k$ and $m \in M$, $m' \leq m$. Therefore, we defined a Kripke structure. This Kripke structure is constrained because any propositional symbol forced at any $k \in K$, $n \leq k$ is also forced at n and every propositional symbols forced at any $m \in M$, $n \leq m$ is forced at n as well.

Neither formula F nor formula G is forced at n in the glued constrained Kripke structure, and hence, formula $F \lor G$ is not forced at n either. Therefore our assumption is incorrect; either F or G should be derivable. ♦

**Theorem.** If formula $F \lor G$ is derivable in LJ+, and every propositional symbol occurring in both

F and G have only positive or only negative occurrences, then either F or G is derivable in LJ+.

**Proof.** The proof of this theorem is by induction on the number of propositional symbols occurring in both F and G. Base case: F and G have no common propositional symbols. The previous lemma gives a proof of the base case. Induction step: suppose this theorem holds for formulas having not more than n propositional symbols occurring in both disjuncts, and F\/G has n+1 propositional symbols occurring in both F and G. Suppose R is one of such propositional symbols, all occurrences of R are positive. The case of negative occurrences is similar to this one. If F\/G is derivable in LJ+, so is F\/G{R|⊥}. Note that F\/G{R|⊥} is the same as F{R|⊥}\/G{R|⊥}. By the induction hypothesis, either F{R|⊥} or G{R|⊥} is derivable in LJ+ because F{R|⊥}\/G{R|⊥} contains n propositional symbols common for F and G.

For certainty, suppose F{R|⊥} is derivable in LJ+. As shown earlier in the section 'Intuitionistic Features', F{R|⊥} ⊃ F is derivable in LJ if all occurrences of R are positive. Due to closeness of LJ+ under modus ponens, F is also derivable in LJ+. ♦

We show that the existence property does not hold in LJ+ by constructing a counterexample. Let A be a predicate symbol and let N be a propositional symbol. Formula ∃x (A(a)\/N)&(A(b)\/¬N)) ⊃A(x) is derivable in LJ+. Here is how it is derived.

```
N |-- N                   A(a) |-- A(a)
========                  =========
N, ¬N |-- A(a)            A(a), ¬N |-- A(a)
-----------------------------------------------
          A(a)\/N, ¬N |-- A(a)
```

Formula A(b)\/¬N, N |-- A(b) is derived in a similar manner. Now, continue the above derivation.

```
(A(a)\/N)&(A(b)\/¬N), ¬N |-- A(a)
-----------------------------------------
¬N |- ((A(a)\/N)&(A(b)\/¬N))⊃A(a)
-----------------------------------------
¬N |- ∃x ((A(a)\/N)&(A(b)\/¬N))⊃A(x)
```

Similarly, one can derive N |- ∃x ((A(a)\/N)&(A(b)\/¬N))⊃A(x) from A(b)\/¬N, N |-- A(b). Finally, application of the neutralization rule gives the formula in question.

In order to show that the existence property does not hold for this formula, we have to consider all possibilities for the variable replacing x: it is a; it is b; it is neither a nor b. Thus, it is sufficient to show that none of the three formulas below is derivable in LJ+:
(A(a)\/N)&(A(b)\/¬N)) ⊃A(a)
(A(a)\/N)&(A(b)\/¬N)) ⊃A(b)
(A(a)\/N)&(A(b)\/¬N)) ⊃A(c)
None of these is even classically derivable. It can be shown by building classical logic interpretations in which these formulas are false. Let the domain of these interpretations be { 0, 1, 2 }. In the interpretation applied to the first and third formula, the set of true atomic values is comprised of N and A(0). Only A(0) is true in the interpretation applied to the second formula. The first formula evaluates to false when a=1, b = 0. The second formula evaluates to false when a=0, b=1. The third formula evaluates to false when a=1, b=0, c=2.

**Theorem**. The existence property holds for the following three classes of formulas derivable in LJ+:
1. $\exists x\, (A(x)\&B)$ if $A(x)$ is unipolar
2. $\exists x\, (A(x)\lor B)$ if $A(x)$ is unipolar, occurrences of every propositional symbol from both $A(x)$ and $B$ are all positive or all negative
3. $\exists x\, (B\supset A(x))$ if $A(x)$ is unipolar, $B$ is purely propositional, occurrences of every propositional symbol from both $A(x)$ and $B$ are all positive or all negative

**Proof**. Formula $\exists x\, (A(x)\&B)$ is equivalent to $\exists x\, A(x)\&B$ in intuitionistic logic. Similarly, formula $\exists x\, (A(x)\lor B)$ is equivalent to $\exists x\, A(x) \lor B$. If $\exists x\, A(x)\ \&\ B$ is derivable in LJ+, then both $\exists x\, A(x)$ and $B$ ought to be derivable in LJ+ due to closeness under modus ponens. Since $\exists x\, A(x)$ is derivable in LJ as well, $A(t)$ is intuitionistically derivable for some variable $t$. Therefore, $A(t)\&B$ is derivable. A similar reasoning applies to the second class. If $\exists x\, A(x) \lor B$ is deriavable in LJ+, then either $\exists x\, A(x)$ or $B$ is derivable by the above theorem about the disjunction property. If $B$ is derivable, then $A(t)\lor B$ is derivable for any $t$. If $\exists x\, A(x)$ is derivable, it is intuitionistically derivable as well. Consequently, $A(t)$ is derivable for some variable $t$, and therefrom $A(t)\lor B$ is derivable. As proved earlier, $B\supset A(x)$ is equivalent to $\neg B\lor A(x)$ in LJ+ provided that $B$ is propositional. Therefore, $\exists x\, (\neg B\lor A(x))$ is derivable in LJ+. Now, the reasoning for the second class can be applied to this case too. Derivability of $\neg B\lor A(t)$ implies derivability of $B\supset A(t)$. ♦

Interestingly, LJ+ has a property that can be classified as a weaker form of the existence property. Note that the term from the statement of the intuitionistic existence property can be extracted from the existence proof. This weaker property has a similar constructiveness flavor – the formula without the existential quantifier can be built algorithmically from the existence proof. Note that in a formulation including functions, word variable should be replaced by word term in the following theorem.

**Theorem**. If $\exists x\, A(x)$ is derivable in LJ+, then for some variables $t_1,\ldots,t_n$, $A(t_1)\lor\ldots\lor A(t_n)$ is derivable in LJ+ as well.

**Proof.** Consider proof of sequent $|\!-\!- \exists x\, A(x)$ without cut. It could be the lower sequent of only two rules: the $\exists$–succedent rule and the neutralization rule. If it is the $\exists$–succedent rule, then the succedent of the upper sequent of this rule gives the formula in question. If it is the neutralization rule, then the upper sequents have the form $B_1,\ldots,B_n\ |\!-\!-\ \exists x\, A(x)$ where $B_1,\ldots,B_n$ are propositional formulas. Let us look at the sequents of this form and such that all sequents below them are also of the same form and such that they are lower sequents of rules in which at least one upper sequent has a different form and its succedent is not a propositional formula. We call them corpus sequents.

Only the two following rules can produce corpus sequents as their lower sequents: thinning in the succedent; $\exists$–succedent. Suppose $A(t_1),\ldots,A(t_n)$ are the upper succedents of the $\exists$–succedent rules resulting in the corpus sequents. There should be at least one $\exists$–succedent rule because otherwise we could derive sequent $|\!-\!-$ . Let us replace the lower succedents of the thinnings in the succedent generating corpus sequents by $A(t_1)\lor\ldots\lor A(t_n)$. Let us replace each of the $\exists$–succedent rules generating corpus sequents by $n-1$ $\lor$-succedent rules resulting in $A(t_1)\lor\ldots\lor A(t_n)$ as the succedent at the bottom. Finally, let us replace all succedents below the corpus ones by $A(t_1)\lor\ldots\lor A(t_n)$.

These transformations produce a correct derivation with endsequent $|\!-\!- A(t_1)\lor\ldots\lor A(t_n)$. Note that only antecedent rules applicable to propositional formuas and the neutralization rule can appear below any of the corpus sequents in the original derivation. Out of these rules, all rules with one upper sequent clearly remain correct after the transformations. As for the rules with two upper sequents, any of them with both upper sequents either lying under corpus sequents or being corpus sequents remains correct because the two upper succedents remain identical. If one upper sequent is not a corpus sequent and not lying under one, then it can only be the ⊃–antecedent rule whose succedent is a propositional formula. And these rules remain correct after the transformations as well. ♦

## Acknowledgement

I would like to thank Sergei Artemov, Grigori Mints and Vladik Kreinovich for their help.